\documentclass[a4paper,12pt,twoside]{article}
\usepackage{amssymb}
\usepackage[final]{lpretex}
\usepackage{title}
\usepackage[title]{appendix}
\usepackage{a4}
\usepackage{amscd}
\input xy
\xyoption{all}

\EndOfDump

\def\DHpartformat#1{(#1) }
\def\DHrefpart#1{(\DHRefpart{#1})}
\LBsetstyleof{partno}{arabic}

\DocVersion[Changed to LaTeX]{1}

\let\define\def

 \def\B {{\mathbb B}} \def\C {{\mathbb C}}
  \def\F {{\mathbb F}}

\def\GG {{\mathbb G}}   
  
  \def\P {{\mathbb P}} 
\def\Q {{\mathbb Q}} \def\R {{\mathbb R}}
\def\T {{\mathbb T}} 
  \def\X {{\mathbb X}}
\def\Z {{\mathbb Z}} 

\define \n {\mathbb N}
\define \z {\mathbb Z}
\define \q {\mathbb Q}
\define \PP {\mathbb P}

\def\sA {{\Cal A}}  \def\sC {{\Cal C}}
 \def\sE {{\Cal E}} \def\sF {{\Cal F}}
  \def\sI {{\Cal I}}
  \def\sL {{\Cal L}}
 \def\sN {{\Cal N}} \def\sO {{\Cal O}}
 
 \def\sT {{\Cal T}} \def\sU {{\Cal U}}
  \def\sX {{\Cal X}}
\def\sY {{\Cal Y}}
\def\sZ {{\Cal Z}}

\def\cZ{{\Cal Z}}
\define \cN {\Cal N}
\define \cf {\Cal F}
\define \cg {\Cal G}
\define \cE {\Cal E}
\define \ce {\Cal E}
\define \cc {\Cal C}
\define \cV {\Cal V}
\define \cA {\Cal A}
\define \cK {\Cal K}
\define \cO {\Cal O}
\define \cF {\Cal F}
\define \cn {\Cal N}
\define \cI {\Cal I}
\define \sP {\Cal P}
\define \sW {\Cal W}
\def\tA {\widetilde{\Cal A}}

\def\bsX{\overline{\sX}}
\def\bsU{\overline{\sU}}
\def\bsZ{\overline{\sZ}}

\def\s {\sigma}

\define \x {\xi}
\define \y {\eta}
\define \G {\Gamma}
\define \r {\rho}
\define \w {\omega}

\def \tU {\widetilde U}

\def\tX {\widetilde X}

\def \tV {\widetilde V}
\def \trho {\widetilde {\rho}}

\def \tp {\widetilde{\mathbb P}}
\define \tH {\widetilde H}
\define \tG {\widetilde{\Gamma}}
\define \tW {\widetilde W}
\define \tF {\widetilde F}
\define \tm {\widetilde m}
\define \St {\widetilde S}
\define \Xt {\widetilde X}
\define \tS {\widetilde S}
\define \tpsi {\widetilde \psi}
\define \txi {\widetilde \xi}
\define \tL {\widetilde L}
\define \tE {\widetilde E}
\define \tl {\widetilde l}
\define \tA {\widetilde A}
\define \tom {\widetilde\omega}
\define \tT {\widetilde T}
\define \tB {\widetilde B}
\define \tf {\widetilde f}
\define \tsA {\widetilde{\sA}}

\define \tM {\widetilde M}
\define \tphi {\widetilde{\phi}}
\define \trho {\widetilde{\rho}}
\define \tR {\widetilde R}
\define \tp {\widetilde p}
\define \tq {\widetilde q}
\define \tc {\widetilde c}
\define \tsF {\widetilde {\sF}}
\define \tsN {\widetilde {\sN}}
\define \tsU {\widetilde {\sU}}
\define \tsX {\widetilde {\sX}}
\define \th {\widetilde h}

\def\pd {\partial}

\def \Dx1 {\frac{\pd}{{\pd} x_1}}
\def \Dy1 {\frac{\pd}{{\pd} y_1}}
\def \Dz1 {\frac{\pd}{{\pd} z_1}}
\def \Dx2 {\frac{\pd}{{\pd} x_2}}
\def \Dy2 {\frac{\pd}{{\pd} y_2}}
\def \Dz2 {\frac{\pd}{{\pd} z_2}}

\def\q {\quad} 

\def\mapdiagr#1{\Big\searrow\rlap{$\raise 5pt\vbox{{\hbox{$\mkern -15mu\scriptstyle#1$}}}$}}   

\def\mapdiagl#1{\llap{$\raise 5pt\vbox{{\hbox{$\scriptstyle#1\mkern
-15mu$}}}$}\Big\swarrow}              

\def\Mapdiagr#1{\nearrow\rlap{$\lower 5pt\vbox{{\hbox{$\mkern
-15mu\scriptstyle#1$}}}$}} 

\def\Mapdiagl#1{\llap{$\lower 5pt\vbox{{\hbox{$\scriptstyle#1\mkern
-15mu$}}}$}\searrow} 

\def\Mapswr#1{\swarrow\rlap{$\lower 5pt\vbox{{\hbox{$\mkern
-15mu\scriptstyle#1$}}}$}}              

\def\Mapnwl#1{\nwarrow\rlap{$\lower 5pt\vbox{{\hbox{$\mkern
-15mu\scriptstyle#1$}}}$}}

\def \inj {\hookrightarrow}

\define \Rhook {\hookrightarrow}

\def \half {\raise1pt\hbox{$\scriptstyle
        \frac{1}{2}\displaystyle$}}

\def \x{{\sl X}\llap{$\mkern -2mu {\scriptstyle -}$}}
\def \End {\operatorname{End}}

\def \Exc {\operatorname{Exc}}

\def \Symm {\operatorname{Sym}}

\def \Pic {\operatorname{Pic}}

\define \Kod {\operatorname{Kod}}
\define \dimension {\operatorname{dim}}
\define \codim {\operatorname{codim}}
\define \contr {\operatorname{contr}}
\define \rk {\operatorname{rank}}
\define \im {\operatorname{im}}
\define \Mor {\operatorname{Mor}}
\define \Cl {\operatorname{Cl}}
\define \Hilb {\operatorname{Hilb}}
\define \degree {\operatorname{deg}}
\define \mult {\operatorname{mult}}
\define \Aut {\operatorname{Aut}}
\define \NS {\operatorname{NS}}
\define \Gal {\operatorname{Gal}}
\define \ch {\operatorname{char}}
\define \Jac {\operatorname{Jac}}
\define \Km {\operatorname{Km}}
\define \Sec {\operatorname{Sec}}
\define \Stab {\operatorname{Stab}}
\define \Br {\operatorname{Br}}
\define \inv {\operatorname{inv}}
\define \tr {\operatorname{tr}}
\define \Frob {\operatorname{Frob}}
\define \Symn {\operatorname{Sym}^n}
\define \Ev {\sE^\vee}
\define \ordp {\operatorname{ord}_p}
\define \Supp {\operatorname{Supp}}
\define \Ann {\operatorname{Ann}}
\define \disc {\operatorname{disc}}
\define \Lie {\operatorname{Lie}}
\define \embdim {\operatorname{embdim}}

\def\bC{\overline{C}}

\def\bsX{\overline{\sX}}
\def\bsU{\overline{\sU}}
\def\sZ{\cZ}
\def\bW{\overline{W}}

\def\hod#1#2#3#4{\ensuremath{\if#30 H^{#2}({#1},{\cal O}_{#1}) \else 
 H^{#2}(#1,\Omega^{#3}\if\relax{#4}\relax_{#1}\else _{#1/#4}\fi)\fi}}

\begin{document}
\title[An exceptional locus]{An exceptional locus in 
the perfect  compactification of 
$A_{\MakeLowercase {g}}$}

\author{N. I. Shepherd-Barron}
\address
{Math. Dept.\\
King's College London\\
Strand\\
WC2R 2LS\\
UK}
\email{N.I.Shepherd-Barron@kcl.ac.uk}
\maketitle
\begin{section}{Introduction}
\medskip
This paper is a sequel to \cite{SB}. The main result of that was a description of
the cone ${\overline{NE}}(A_g^P)$ of curves on the perfect compactification
of the coarse
moduli space $A_g$ (the notation will be explained below). In particular, 
the $\Q$-divisor
class $L_g=12H_g-D_g$ is nef but not ample. 
Here we make this result more explicit in two ways:
we identify the exceptional locus $\Exc(L_g)$
as the locus $A_{1,g-1}^P$ 
inside $A_{g}^P$ that is the closure of the locus in $A_g$ that
parametrizes abelian varieties with an elliptic factor,
and we show that $L_g$ is semi-ample in positive characteristic.
We do not know whether $L_g$ is semi-ample in characteristic zero.

Recall that $A_g$ is known to be of general type if $g\ge 7$
and the characteristic is zero, and that 
the canonical class is given by $K_{A_g^P}\sim (g+1)H_g-D_g$. 
It follows (Corollary \ref{interior})
that, if $7\le g\le 10$, then $K_{A_g^P}$ is not nef and therefore the first step in 
running the characteristic zero Minimal Model Program
on $A_g^P$ is a flipping contraction that, in particular, 
crushes the locus that parametrizes
products of $g$ elliptic curves to a point. In particular, each non-trivial
fibre of this contraction penetrates $A_g$, which we regard as
the interior of $A_g^P$. Therefore to run the MMP
on $A_g^P$ (which has terminal singularities if $g\ge 6$;
see [AS] for a proof of this that completes the argument in \cite{SB})
while maintaining a modular interpretation would necessitate
changing the definition of ``principally polarized abelian variety''
in the presence of elliptic factors. In contrast, there is
a prospect that 
the MMP for the coarse moduli space $\overline{M}_g$ 
of stable curves of genus $g$ will run in a way
that does not touch the interior of $\overline{M}_g$,
where in this context ``interior'' refers to the locus $M_g$
of smooth curves [HH]. 

Here is some of the notation used in this paper.
We let $\sA_g$ denote the stack, and $A_g=[\sA_g]$ the coarse
moduli space,
of principally polarized abelian $g$-folds
and $\sA_g^P$, $A_g^P$ their perfect compactifications.
These are particular toroidal compactifications, and dominate the
Satake compactification
$A_g^{Sat}$ of $A_g$.
Stacks that classify abelian varieties with level $n$ structures
will be denoted by $\sA_{g;n}$, etc. In general the geometric quotient
of an algebraic stack $\sX$ will, when it exists, be denoted by $[\sX]$.

We let $H_g$
denote the line bundle of weight $1$ modular forms on $\sA_g^P$.
It gives
a $\Q$-line bundle on $A_g^P$ that can, in turn, be identified
with the pullback of a $\Q$-line bundle on $A_g^{Sat}$.
The boundary $D_g=\sA_g^P\setminus\sA_g$ is a geometrically
irreducible $\Q$-Cartier divisor, but is not Cartier;
even at the stack level there are denominators.


The objects $\sA_g^P$ and $A_g^P$ exist over $\Sp\Z$ and there 
is a contraction $\pi:\sA_g^P\to A_g^{Sat}$
that factors through $A_g^P$; see \cite{FC} for this. As mentioned above,
the main result of
\cite{SB} is a description of the cone of curves ${\overline{NE}}$
(in Mori's sense) of $A_g^P$ over any field:
${\overline{NE}}(A_g^P)$ is the rational cone spanned by certain curves
$C_1$ and $C_2$. Here $C_1$ is the closure of the locus of
products $E\times B$, where $E$ is a varying elliptic curve and
$B\in A_{g-1}$ is fixed (so that $C_1$ is a copy of $A_1^P$,
the compactified $j$-line) and $C_2$ is any curve in
the boundary divisor $D_g=A_g^P\setminus A_g$ such that 
$\pi(C_2)$ is a point.

In consequence \cite{SB}, 
$aH_g-D_g$ is ample on $A_g^P$ (or on $\sA_g^P$)
if and only if $a>12$, and is nef
if and only if $a\ge 12$. 
(I am grateful to Alexeev for pointing
out that the notion of an ample invertible sheaf 
makes sense on any proper stack $\sX$ with finite inertia
and geometric quotient $X=[\sX]$:
the invertible sheaf $\sL$ on $\sX$ is ample if and only if,
for all coherent sheaves $\sF$ on $\sX$,
$H^i(\sX,\sL^{\otimes n}\otimes\sF)=0$ for all $i>0$ and
for all $n\gg 0$. If $\sL$ corresponds to the $\Q$-line bundle $L$ on $X$,
then $\sL$ is ample if and only if $L$ is so. Note also that, even more obviously,
the notions of nef and semi-ample are insensitive to the distinction
between $\sX$ and $X$, and $\Exc(\sL)$ is exactly the inverse image in $\sX$
of $\Exc(L)$.)

Recall \cite{Ke} that if $L$ is a nef divisor class
on a projective variety $X$, then the \emph{exceptional locus} 
$\Exc(L)$ of $L$ is the union of all the subvarieties
$Z$ of $X$ with $L^{\dim Z}.Z=0$, and 
that $L$ is semi-ample on $X$
provided that it is semi-ample on $\Exc(L)$
and the ground field $k$
has $\ch k>0$. Notice that Keel's result holds without change
when $X$ is allowed to be any proper stack over a field
whose inertia stack is finite (so that a geometric quotient of $X$ exists
as a proper algebraic space, by \cite{KM})
and whose geometric quotient is projective.

Here is the main result of this paper; I am very grateful to
Stefan Schr{\"o}er, who asked whether $12H_g-D_g$ is semi-ample
in positive characteristic.

\begin{theorem}\label{exceptional} 
(= Theorem \ref{4.4})
Suppose that the base is a field $k$
and set $L_g=12H_g-D_g$.

\part [i] $\Exc(L_g)$ is the subvariety $A_{1,g-1}^P$
that is the image of the stack
$\sA^P_{1,g-1}$ defined in Corollary (\ref{defn of product}).

\part [ii] Suppose \emph{either}
that $\ch k=0$ and that $g\le 11$
\emph{or} that $\ch k>0$ and that $g$ is arbitrary.
Then $L_g$ is semi-ample.
For a sufficiently divisible integer $n$ the linear system
$\vert mL_g\vert$
defines a birational morphism 
$\phi_{\vert mL_g\vert}=\contr_{R_1}:A_g^P\to V_g$
onto a normal projective
$k$-variety $V_g$ that contracts the ray
$R_1$ generated by the curve $C_1$.

\part[iii] The normalization of the
image $\phi_{\vert mL_g\vert}(\Exc(L_g))$ 
is isomorphic to $V_{g-1}$.
\noproof
\end{theorem}
Observe that
$V_g$ can be regarded as being obtained from $A_g^P$ by 
``crushing elliptic factors to points''. From this point of view
it is not clear \emph{a priori}
that $V_g$ exists as an algebraic variety; indeed,
we are unable to prove such a statement when $g\ge 12$
and the characteristic is zero.
\end{section}
\begin{section}{A remark and an example}
\medskip
Given a principally polarized abelian variety $(A,\lambda)$
there is not usually a theta divisor $\theta$ on $A$,
although there is a natural line bundle $\sL_{2\lambda}$
associated to $2\lambda$. However, the stack $\sA_g$
is naturally isomorphic to the stack $\sX_g$
of pairs $(X,\theta)$ where $X$ is a symmetric torsor under $A$ 
and $\theta$ is an ample divisor on $X$ that
defines $\lambda$ on $A=\Aut^0_X$. 
Moreover, $(X,\theta)$ has a $1$-dimensional factor
if and only if $(A,\lambda)$ has an elliptic factor.

It is sometimes convenient
to identify $\sX_g$ with $\sA_g$
and to consider pairs $(X,\theta)$
instead of ppavs $(A,\lambda)$.

Now we recall \cite{{H},{DO}} some of the classical geometry of the moduli
space $A_2$ over $\C$.

A level $2$ structure on a principally polarized abelian surface
$(A,\lambda)$ is a symplectic isomorphism 
$\psi:A[2]\to G:=(\Z/2)^2\times\mu_2^2$,
where $G$ has its standard symplectic pairing
and $A[2]$ has the Weil pairing defined by the principal polarization
$\lambda$. There is a standard projective action of $G$ on $\P^3$;
the linear system $\vert 2\theta\vert$ on $X$
then gives a $G$-equivariant morphism $X\to\P^3$ that factors through the Kummer surface
$\Km(X)=X/(-1)$. 
The image of $\Km(X)$ lies in a unique $G$-invariant
quartic; taking the coefficients of this quartic then determines
a point on $\Sigma$, the Segre cubic threefold \cite{H}; this is the unique
cubic threefold with $10$ nodes and lies in $\P^5$ with equations
$e_1=e_3=0$, where $e_i$ is the $i$th elementary symmetric function
in $6$ variables. If $(X,\theta)$ is irreducible
then $\vert2\theta\vert$ is very ample on $\Km(X)$, so $(X,\theta)$
is determined by the point on $\Sigma$, while if $X=E_1\times E_2$,
where $E_1,E_2$ are curves of genus $1$,
then the image of $\Km(X)$ is $\P^1\times \P^1$, which does not determine
$X$.
\def\tSigma{\widetilde\Sigma}

There are also $15$ planes on $\Sigma$. Let $\tSigma\to \Sigma$ be the 
blow-up of the nodes, $\Pi=\sum \Pi_i$ the exceptional divisor of the blow-up
and $D=\sum D_i$ the strict transform of the sum of the planes.

If $G$ is a finite group acting on a variety $V$ then we let $V/G$
denote the geometric quotient.

\begin{theorem} \part[i]
$\tSigma=A^P_{2;2}$, the perfect compactification of the
level $2$ moduli space $A_{2;2}$. $D$ is the toroidal boundary and
$\Pi$ the locus of products $E_1\times E_2$ of elliptic curves. 

\part[ii] $A^P_2$ has two distinct birational contractions: one
is the standard contraction $\pi:A^P_2\to A_2^{Sat}$ and the
other is a birational contraction $\rho:A_2^P\to \Sigma/\frak S_6$, where
$\frak S_6$ is the symmetric group on $6$ letters and is isomorphic to
$Sp_4(\F_2)$. 

\part[iii] $\Sigma/\frak S_6$ is isomorphic to a weighted projective space
$\P(2,4,5,6)$; $A_{2;2}^{Sat}$ is isomorphic to the Igusa quartic
$e_1=e_4-e_2^2=0$ in $\P^5$ (this is also the projective dual of $\Sigma$);
and $A_2^{Sat}$ is isomorphic to $\P(2,3,5,6)$.
\noproof
\end{theorem}
From one point of view, 
the point of Theorem \ref{exceptional} 
is to extend this picture, albeit in a more abstract way, to all
values of $g$.
\end{section}
\begin{section}{The structure of the perfect boundary}
For any $g$ let $\Lambda_g$ be a fixed copy of $\Z^g$
and $B(\Lambda_g)$ the $\Z$-module of symmetric
bilinear $\Z$-valued bilinear forms on $\Lambda_g$.
Let $C_g$ denote the cone
of real positive definite symmetric bilinear forms
in $g$ variables and $\bC_g$ the cone of
real positive semi-definite symmetric bilinear forms
in $g$ variables. We identify both of these
as subsets of $B(\Lambda_g)\otimes_\Z\R$. Identify $C_g$
with the interior of $\bC_g^0$ of $\bC_g$. A toroidal
compactification of $\sA_g$ or $A_g$ corresponds to a choice
of \emph{$GL_g(\Z)$-admissible} decomposition of $\bC_g$, as defined in
\cite{AMRT}. We consider here
a particular admissible decomposition $\Sigma^P$ of $\bC_g$, namely, that defined
by the perfect quadratic forms. It is called the {\emph{perfect decomposition}}
and the cones appearing in it (that is, the facets) are the
\emph{perfect cones}. It is described in Theorem \ref{voronoi}
below.
Since our results are particular to this specific toroidal
compactification,
the arguments of this paper
depend upon knowing something about
perfect forms and the combinatorics of the
decomposition of $\bC_g$ that they define.
The cones of maximal dimension $g(g+1)/2$
correspond to perfect quadratic forms
(this can be taken as the definition of a perfect quadratic form).

\begin{theorem}\label{combinatorics}
\part[i] (Voronoi) \cite{AMRT} \label{voronoi}
The convex hull of the positive semi-definite integral
forms in $\bC_g$ defines a $GL_g(\Z)$-admissible decomposition
$\Sigma_g^P$ of $\bC_g$.


\part[ii] (Barnes-Cohn) \cite{BC}\label{bc} This convex hull coincides with
the convex hull of the primitive rank $1$ forms. Moreover,
a form of rank at least $2$ lies in the interior of the hull.

\part[iii] \label{product} If $p=p(x_1,...,x_m)$ and $q=q(y_1,...,y_n)$
are perfect forms of equal minimal norm, then
$p+q=r=r(x_1,...,x_m,y_1,...,y_n)$ is a form that defines
a perfect cone in
$\bC_{m+n}$.

\part[iv] \label{closure} Suppose that $\tau$ is a perfect cone in $\bC_r$
that meets $C_r$.
Then its closure contains perfect cones in some copy of $\bC_{r-1}$.
In particular, if $\tau$ is a minimal perfect cone in $\bC_r$ that
meets $C_r$, then $\tau\cap\bC_{r-1}$ is a union of minimal
perfect cones in $\bC_{r-1}$ that meet $C_{r-1}$.

\part[v] \label{codimension} Suppose that $\sigma$ is a maximal perfect cone
in $\bC_{r}$, so that 
$\sigma$ is defined by
a perfect form $q$ in $r$ variables. Suppose also that $\tau_1,\tau_2$
are closed cones in $\bC_{r+1}$ such that both contain
$\sigma$ in their boundary and $\dim\tau_i=\dim\sigma+1$.
That is, both cones $\tau_i$ are minimal with respect to 
containing $\sigma$ in their boundary. Then each
$\tau_i$ can be defined by a quadratic form $\lambda q+l_i^2$,
where $l_i$ is a primitive linear form and $\lambda\in \Q$
is the inverse of the minimal norm of $q$.
In particular, $\tau_1$ and $\tau_2$ are conjugate
under the parabolic
subgroup of $GL_{r+1}(\Z)$ that preserves each of the
first $r$ variables.

\begin{proof} For \DHrefpart{i} and \DHrefpart{ii} see the references.
\DHrefpart{iii} is trivial. For \DHrefpart{iv}, we can suppose
that $\sigma$ is minimal with respect to meeting $C_r$. Suppose
that $l_1,...,l_n$ are primitive elements of $\Lambda_r^\vee$
such that $l_1^2,...,l_n^2\in B(\Lambda_r)$ span the $1$-dimensional
faces of $\sigma$. So there are $\lambda_1,...,\lambda_n>0$ such that
$\sum \lambda_il_i^2\in C_r$. Then, for some $t$, 
$\sum_1^t \lambda_il_i^2$ is of rank at most $r-1$ for every 
$\lambda_1,...,\lambda_t>0$, while
$\sum_1^{t+1}\lambda_il_i^2$ is of rank $r$ for some 
$\lambda_1,...,\lambda_{t=1}>0$.
Take $t$ to be maximal subject to this; then the
cone $\tau$ generated by $l_1^2,...,l_t^2$ is a face of $\sigma$
with the stated properties.

\DHrefpart{v} 
is even easier.
\end{proof}
\end{theorem}

We denote by $\sA_g^P$ the toroidal compactification of $\sA_g$
that corresponds, as in [FC], to the perfect decomposition
of $\bC_g$.
The boundary $D=\sA_g^P\setminus\sA_g$ 
is a divisor, and is the inverse image $\pi^{-1}(A_{g-1}^{Sat})$,
where $A_{g-1}^{Sat}$ is identified with the boundary
$A_g^{Sat}\setminus A_g$.
We also have the partial compactification $\sA_g^{part}$, which is open
in $\sA_g^{tor}$ and is, by definition, the inverse image $\pi^{-1}(A_g\coprod A_{g-1})$.
The universal abelian
scheme $\sU_{g-r}\to\sA_{g-r}$ has an extension to a semi-abelian
scheme $\sU_{g-r}^{part}\to \sA_{g-r}^{part}$ whose degenerate
fibres have torus rank $1$. Taking $r$-fold fibre products
gives a semi-abelian scheme $\delta:\sU_{g-r}^{r,part}\to
\sA_{g-r}^{part}$

In higher codimension the boundary is described as follows.
There is a stratified scheme $F=F_r$, locally of finite type over the base,
the closures of whose strata are projective toric varieties,
with an action  of $GL_r(\Z)$ on $F_r$, and an $F_r$-bundle
$\sF_r\to\sU_{g-r}^r$, with an equivariant action of $GL_r(\Z)$,
such that $\pi^{-1}(A_{g-r})=\sF_r/GL_r(\Z)$, the stack quotient.
The closures $F_{r,\tau}$ of the various strata of $F_r$
correspond to the perfect cones $\tau$ in the perfect decomposition of
$\bC_r$ that meet
the interior $C_r$. 
In particular, the irreducible components of $F_r$ correspond to the
minimal such cones.

Here is a more precise description of $F_r$: put $\Lambda_r=\Z^{\oplus r}$,
let $B_r$ be the lattice of symmetric bilinear forms on $\Lambda_r$
and denote by $T_r$ the torus with character group
$\X^*(T_r)=B_r$. Then there
is a locally finite $GL_r(\Z)$-equivariant torus embedding $T_r\inj Y_r$ 
such that $F_r$ is a $T_r\rtimes GL_r(\Z)$-equivariant
closed subscheme of the boundary $Y_r\setminus T_r$. 
The closure
$F_{r,\tau}$ of a stratum in $F_r$ is then a torus embedding under a quotient $T$
of $T_r$
and gives rise to the closure $\sF_{r,\tau}$ of a stratum in
$\sF_r$; this closure is a proper $F_{r,\tau}$-bundle $\sF_{r,\tau}\to\sU_{g-r}^r$
that is a relative $T$-equivariant compactification of a $T$-bundle
$\sT\to\sU_{g-r}^r$.

In turn, the image
of $\sF_{r,\tau}$ in $\sA_g^P$ is an irreducible closed substack
$\sX_{r,\tau}$ of
$\pi^{-1}(A_{g-r})$. If $n\ge 3$ and
is invertible in the base, then in the stack $\sA_{g;n}^P$ the
image $\sX_{r,\tau}$ can be identified
with $\sF_{r,\tau}$. In this case (that is, at level $n$)
$T_r$ acts on $\sX_{r,\tau}$, via the quotient $T_r\to T$, 
and this action extends to an action on the closure $\bsX_{r,\tau}$ 
of $\sX_{r,\tau}$ in $\sA_{g;n}^P$.

Each such $\tau$ lies in the closure of finitely many maximal perfect
cones $\s$ in $\bC_r$. Such a cone $\s$ corresponds to the choice, up to scalars, of a perfect form $q$
in $r$ variables.
The closure
$\bsU^r_{g-r,\s}$ of $\sU^r_{g-r}$ is just
$\bsX_{r,\s}$.
We let $\bsU_{g-r,\s}^{r,norm}$ denote the normalization
of $\bsU^r_{g-r,\s}$.

Translating Theorem \ref{combinatorics} into algebraic geometry
yields the following statement.

\begin{corollary}\label{alg_geom}
\part[i] 
(Corollary of \ref{voronoi}) There exist
toroidal compactifications $\sA_g^P$ and $A_g^P$
corresponding to this decomposition.

\part[ii] (Corollary of
\ref{bc}) As a Deligne-Mumford stack, $\sA_g^P$ has terminal singularities
and the boundary $D$ is absolutely irreducible.

\part[iii]\label{defn of product}
(Corollary of \ref{product}) The product morphism
$\sA_g\times \sA_h\to \sA_{g+h}$
extends to a morphism $\sA_g^P\times \sA_h^P\to \sA_{g+h}^P$
whose image is denoted by $\sA^P_{g,h}$.

\part[iv] (Corollary of \ref{closure}) 
$\pi^{-1}(A_{g-r}^{Sat})$ lies in the closure of
$\pi^{-1}(A_{g-r+1})$ and $\pi^{-1}(A_{g-r+1}^{Sat})$ is the closure
of
$\pi^{-1}(A_{g-r+1})$.

\part[v] (Corollary of \ref{codimension}) A maximal perfect cone $\sigma$ in
$\bC_r$ corresponds to an irreducible closed substack $\bsU_{g-r,\s}^r$ of
$\pi^{-1}(A_{g-r}^{Sat})$ that contains
an open substack isomorphic to $\sU_{g-r}^{r, part}$. 
The complement
$\bsU_{g-r,\s}^r\setminus\sU_{g-r}^{r, part}$ has codimension 
at least $2$ in $\bsU_{g-r,\s}^r$.
\begin{proof}
\DHrefpart{i}, \DHrefpart{ii} and \DHrefpart{iii} are immediate.

For \DHrefpart{iv}, we use the fact that the irreducible components $\sZ$ of 
$\pi^{-1}(A_{g-r})$ correspond to the (equivalence classes of the) minimal
perfect cones $\tau$ in $\bC_r$ that meet $C_r$. Since, by \ref{closure},
$\tau\cap\partial\bC_r$ is a union of minimal perfect cones in
$\bC_{r-1}$ that meet $C_{r-1}$, the closure $\bsZ$ of $\sZ$
lies in the closure of an irreducible component of
$\pi^{-1}(A_{g-r+1})$.
So $\pi^{-1}(A_{g-r})\subseteq {\overline{\pi^{-1}(A_{g-r+1})}}$. Then,
by induction, $\pi^{-1}(A_{g-r-m})\subseteq
{\overline{\pi^{-1}(A_{g-r+1})}}$,
so that $\pi^{-1}(A_{g-r}^{Sat}) \subseteq {\overline{\pi^{-1}(A_{g-r+1})}}$.
Therefore $\pi^{-1}(A_{g-r+1}) =  {\overline{\pi^{-1}(A_{g-r+1})}}$.

For \DHrefpart{v}, note that the irreducible components of 
$\bsU_{g-r,\s}^r\setminus\sU_{g-r}^{r, part}$ correspond to the minimal
perfect cones $\tau$ in $\bC_{r+1}$ that contain $\s$. These are
equivalent, and we are done.
\end{proof}
\end{corollary}

\begin{definition} An open substack $\sU$ of an algebraic stack $\sX$
is \emph{nearly equal} 
(abbreviated to
n.e.) to $\sX$ if its complement has codimension
at least $2$ everywhere. A \emph{near equality}
is an open embedding $\sU\inj\sX$ of stacks
whose image is nearly equal to its target. 
\end{definition}

In particular, $\sA_g^{part}$ is nearly equal to $\sA_g^P$.
\end{section}
\begin{section}{The exceptional locus of $12H_g-D_g$}\label{exc_locus}
Set $L_g=12H_g-D_g$.

\begin{proposition}\label{no base points}
For sufficiently divisible $m$,
the linear system $\vert mL_g\vert$ has no base points in
$\sA_g^{part}$ and contracts all curves (that is, complete 
$1$-dimensional substacks) of the form
$\sA_1^P\times\{B\}$ where $B$ is a point in $\sA_{g-1}$. 
This morphism separates points
except along $A_{1,g-1}^P\cap A_g^{part}$.
\begin{proof}
First, work over $\Z[1/2]$ and impose a level $2$ structure.
We then consider the morphism defined by the $2\theta$
linear system: there is a universal family
$f:V\to\sA_{g;2}^{part}$ of projective schemes
with level $2$ structure, and the $2\theta$
linear system defines, by taking the cycle-theoretic
image of each fibre of $f$, a morphism $\Phi$ from
$\sA_{g;2}^{part}$ to the Chow scheme of
$\P^{2^g-1}$. After dividing by the finite group
$Sp_{2g}(\Z/2)$ we then get a morphism
$\phi$ from $\sA_g^{part}$ to some scheme,
and $\phi$ contracts every curve of the form
$A_1^P\times\{B\}$, since the Kummer variety of an elliptic
curve $E$ is independent of $E$, so the $j$-line collapses to a point.
It follows that $\phi$ is defined by some linear system
$\vert mL_g\vert$, since when $g=1$ the bundle $H_g$ has degree 
$1/12$ and $D_g$ has degree $1$.

Over $\Z[1/3]$ we work at level $3$.
There is a normalized projective space $\P^{3^g-1}$ 
such that, given $(X,\theta)$, $X$ is embedded
in $\P^{3^g-1}$ via $\vert 3\theta\vert$.
Let $Gr$ denote the Grassmannian of quadrics in $\P^{3^g-1}$.
Then to $(X,\theta)$ associate the point $P_{(X,\theta)}$ 
in $Gr$ that
corresponds to the space of quadrics that pass through the
image of $X$ in $\P^{3^g-1}$.
Since curve of genus $1$ and degree $3$ cannot be recovered
from the quadrics that contain it (there are no such quadrics), this 
association defines a morphism $\sA_{g;3}^{part}$ that performs
a similar function of ``losing elliptic factors''.
\end{proof}
\end{proposition}

So
pick $r\ge 2$ and consider $L_g$ on the inverse image
$\pi^{-1}(A_{g-r})$, where $A_{g-r}$ is regarded as a stratum
in $A_g^{Sat}$. The closure $\sZ_r$ of $\pi^{-1}(A_{g-r})$ 
in $\sA_{g}^P$ is a finite union of irreducible components $\bsX_{r,\tau}$
as above.
Note that $\sZ_r=\pi^{-1}(A_{g-r}^{Sat})$, as already pointed out.

Fix a maximal perfect cone $\sigma$ in $\bC_r$ and a face $\tau$ of
$\sigma$ that meets $C_r$. There is, up to scalars, a unique
perfect quadratic form $q$ which defines $\s$ and
can be written, in many ways, as a linear combination
$q=\sum_1^r\lambda_ix_i^2$
where $x_i$ is a primitive integral linear form and $\lambda_i\in\Q_+$.
Each rank $1$ form $x_i^2$ corresponds, in 
$\End(\sU_{g-r}^r)\otimes\Q$, to a projection
$\xi_i:\sU_{g-r}^r\to \sU_{g-r}$ over $\sA_{g-r}$.
So (since homomorphisms of abelian schemes extend uniquely
to homomorphisms of semi-abelian schemes) there is a commutative diagram
$$\xymatrix{
{\bsX_{r,\tau}}\ar@{<-^{)}}[r]^{cl.} 
&\ {\bsU^r_{g-r,\s}} \ar@{<-^{)}}[r]^{n.e.} &\  
{\sU_{g-r}^{r,part}}\ar[dll]^{\alpha_{g,r}} \ar[r]^{\xi_i} \ar[d]_{\delta}  &
\ {\sU^{part}_{g-r}}\ \ar@{^{(}->}[r]^{n.e.} \ar[dl]^{\gamma}  & 
{D_{g-r+1}}\ \ar@{^{(}->}[r]^{cl.}& {\sA_{g-r+1}^P.}\\
{\sA_{g}^P}\ \ar@{<-^{)}}[u]^{cl.}&& {\sA_{g-r}^{part}}
}$$
(Here the labels $cl.$ and $n.e.$ refer, respectively,
to closed embeddings and near equalities.)

\begin{proposition}\label{5.3} 
Let $s$ denote the $0$-section of $\delta$. Then the restrictions
$\alpha^*_{g,r}D_g\vert_s$ and $\delta^*D_{g-r}\vert_s$ are linearly
equivalent.
\begin{proof}
Recall that the multiplication $\sA_{g-r}\times\sA_r\to\sA_g$ extends
to $\sA_{g-r}^P\times\sA_r^P\to\sA_g^P$.
The perfect form $q$ corresponds to a maximally degenerate
boundary point $x\in \sA_r^P$, and then the image
of $\{x\}\times \sA_{g-r}^P$ in $\sA_g^P$
is exactly the closure in $\sA_g^P$ of the zero section
of $\sU_{g-r}^r\to\sA_{g-r}$.

Identify $s$ with $\{x\}\times \sA_{g-r}^{part}$
as above. Then there is a chain of rational curves in $\sA_r^P$
leading from $x$ to an interior point $y=E^r$, for any elliptic curve
$E$ (let $E$ degenerate, and then take a chain of rational curves leading from
this maximally degenerate boundary point to $y$).

So $D_g\vert_s$ is rationally equivalent to
$D_g\vert_{\{y\}\times \sA_{g-r}^{part}}$,
and now the result is obvious.
\end{proof}
\end{proposition}

Set $\Lambda=-\alpha^*_{g,r}D_g+\delta^*D_{g-r}$, a divisor
class on $\sU_{g-r}^{r,part}$. Since
$\sU_{g-r}^{r,part}\to\sA_{g-r}^{part}$ is semi-abelian and $\Lambda$
is trivial on the zero section, $\Lambda$
is determined by the polarization that it defines on the generic fibre
$\sU_{g-r,\eta}^r$.
(At this point we also use the fact that line bundles on the
semi-abelian scheme
$\sU_{g-r}^{r,part}\to\sA_{g-r}^{part}$ that are trivial on the zero section are
determined by the polarization that they define on the generic fibre.
This is because the only global point of the self-dual abelian scheme
$\sU_{g-r}^r\to\sA_{g-r}$ 
is zero.)

\begin{corollary}
$$\Lambda\sim \sum_1^r\lambda_i\xi_i^*\left(-D_{g-r+1}\vert_{\sU_{g-r}^{part}}+\gamma^*D_{g-r}\right).$$
\begin{proof}
From its definition, and knowledge of the polarization defined by
$-D_g$ on $\sU_{g-r,\eta}^r$, 
the polarization defined by $\Lambda$ on
$\sU_{g-r,\eta}^r$ is the quadratic form $q$. 
The polarization defined by $-D_{g-r+1}$ on the generic fibre of
$\sU_{g-r}$ is that given by the primitive rank one form $x_i^2$,
and we are done.
\end{proof}
\end{corollary}

For $i=1,\ldots,r$, choose a large positive integer $n_i$ and set
$\txi_i=[n_i]\circ\xi_i:\sU_{g-r}^{r,part}\to\sU_{g-r}^{part}$ and 
$\rho_i = \alpha_{g-r+1,1}\circ\txi_i:\sU_{g-r}^{r,part}\to
\sA_{g-r+1}^P$. Then
$$\alpha_{g,r}^*\Lambda\sim \sum_1^r\frac{\lambda_i}{n_i^2}\txi_i^*\left(-D_{g-r+1}\vert_{\sU_{g-r}^{part}}+\gamma^*D_{g-r}\right),$$
which can be re-written as
$$L_g\vert_{\sU_{g-r}^{r,part}}\sim\left(1-\sum_i\frac{\lambda_i}{n_i^2}\right)\delta^*L_{g-r}+\sum_i\frac{\lambda_i}{n_i^2}
\rho_i^*L_{g-r+1}.$$
We abbreviate this to
\begin{equation}\label{(*)}
L_g\vert_{\sU_{g-r}^{r,part}}\sim a \delta^*L_{g-r} +\sum b_i\rho_i^*L_{g-r+1}.
\end{equation}
Note that $a,b_i>0$.

\begin{theorem}\label{4.4} 
(= Theorem \ref{exceptional}) \part [i] \cite{SB} $L_g$ is nef.

\part [ii] $\Exc(L_g)=\sA^P_{1,g-1}$.

Now suppose also that
\emph{either}
that the characteristic is zero and that $g\le 11$
\emph{or} that the characteristic is positive.

\part [iii] 
Then $L_g$ is semi-ample and 
there is a normal projective $k$-variety $V_g$
and a birational contraction $A_g^P\to V_g$ of
the ray $R_1$ whose exceptional
locus is $A_{1,g-1}^P$.

\part[iv] The normalization of the
image of $A_{1,g-1}^P$ is isomorphic to $V_{g-1}$.
\begin{proof} \DHrefpart{i} was proved in \cite{SB}. However, in order to
prove \DHrefpart{ii} 
it is simplest to give here another proof of \DHrefpart{i}
whose techniques can be extended to prove \DHrefpart{ii} also.

Assume that $L_g$ is not nef, and
that $\sC$ is a complete
curve in $\sA_g^P$ with $L_g.\sC<0$. By Proposition \ref{no base points}, and
since $-D_g$ is $\pi$-ample, there is some
minimal $r\ge 2$ such that $\sC$ maps to a curve in $A_{g-r}^{Sat}$
that does not lie in $A_{g-r-1}^{Sat}$.

Fix an integer $n\ge 3$ that is prime to $\ch k$
and consider the perfect compactification
$\sA_{g;n}^P$ of the level $n$ stack $\sA_{g;n}$.

Consider $L_g$ on the inverse image
$\pi^{-1}(A_{g-r})$, where $A_{g-r}$ is regarded as a locally closed
subvariety
of $A_g^{Sat}$. The closure $\cZ^r$ of $\pi^{-1}(A_{g-r})$ 
in $\sA_{g}^P$ is a finite union of irreducible components
$\bsX_{r,\tau}$, each of which is the image of an equivariant closure
$\bsX_{r,\tau;n}$
in $\sA_{g;n}^P$ of a 
$T$-bundle over the universal level $n$ abelian scheme $\sU_{g-r;n}^r$, 
where $T$ is a quotient 
of the torus $T_r$ whose cocharacter group is $B_r$.

We can use the $T$-action at level $n$ to construct
a specialization (that is, a 
rational equivalence) $\sC\sim \sC_0+\sF$, where $\sC_0$ is contained 
in a minimal stratum of $\cZ_r$ and $\sF$
is the image in $\sA_g^P$
of a closed subvariety $\sF_n$ in $\sA_{g;n}^P$
such that $\sF_n$ is preserved by $T$ but no component of $\sF_n$
consists of $T$-fixed points.
So then $\sF$ is $\pi$-vertical, so that
$(-D_g).\sF>0$ and $H_g.\sF=0$, and then $\sC_0$ is $L_g$-negative.
So we can assume that $\sC$ lies in a minimal stratum of $\cZ_r$.
Recall that each such stratum is one of the
closed substacks $\bsU^r_{g-r,\s}$ considered previously,
and so $\sC\subseteq \bsU^r_{g-r,\s}$, say.

Note that, for any sufficiently divisible integer $m$,
there are equalities
\begin{eqnarray*}
&&H^0(\bsU^{r,norm}_{g-r,\s},\sO(mL_g)\vert_{\bsU^{r,norm}_{g-r,\s}})
=H^0(\sU_{g-r}^{r,part},\sO(mL_g)\vert_{\sU_{g-r}^{r,part}})\\
&=& H^0(\sU_{g-r}^{r,part},\sO(ma\delta^*L_{g-r})\otimes
H^0(\sU_{g-r}^{r,part},\sO(\sum mb_i\rho_i^*L_{g-r+1})).
\end{eqnarray*}
The first equality follows from the facts that
$\sU_{g-r}^{r,part}$ is nearly equal to $\bsU^{r,norm}_{g-r,\s}$
and that $\bsU^{r,norm}_{g-r,\s}$ is normal.
The second follows from the linear equivalence \ref{(*)}
and the fact that, by
Proposition \ref{no base points},
$L_h$ has (stably) no base points on $\sA_h^{part}$
for all values of $h$, and in particular for
$h=g-r$ and $h=g-r+1$.

Therefore, the stable base locus
$\B(L_g\vert_{\bsU^{r,norm}_{g-r,\s}})$ is contained
in the boundary $\bsU^{r,norm}_{g-r,\s}\setminus\sU^{r,part}_{g-r}$
of $\bsU^{r,norm}_{g-r,\s}$.

The next lemma is well known but we include a proof for lack of a 
convenient reference.

\begin{lemma} If $\nu:X\to Y$ is a finite dominant morphism of integral proper
algebraic spaces and $L\in\Pic_Y$, then $\nu^{-1}\B(L)=\B(\nu^*L)$.
\begin{proof} The homomorphism
$$R(Y,L)=\oplus H^0(Y,L^{\otimes n})\to R(X,\nu^*L)=\oplus H^0(X,\nu^*L^{\otimes n})$$
is injective and finite.

Suppose that $x\in\B(\nu^*L)$, so that
$s(x)=0$ for all $s\in R(X,\nu^*L)_{{}\ge 1}$. 
In particular,
$t(x)=0$ for all $t\in R(Y,L)_{{}\ge 1}$
and so $x\in\nu^{-1}\B(L)$. That is,
$\B(\nu^*L)\subseteq\nu^{-1}\B(L)$.

Now suppose that $x\in\nu^{-1}\B(L)\setminus\B(\nu^*L).$
Choose $s\in H^0(X,\nu^*L^{\otimes n})$ with $s(x)\ne 0$.
There is an equation
$$s^N+a_{N-1}s^{N-1}+\cdots+a_0=0$$
for some $a_i\in H^0(Y,L^{\otimes n(N-i)})$.
Then $a_i(x)=0$ for all $i$, which is absurd.
\end{proof}
\end{lemma}

It follows that $\B(L_g\vert_{\bsU^r_{g-r,\s}})$
is contained in 
the boundary $\bsU^{r}_{g-r,\s}\setminus\sU^{r,part}_{g-r}$
of $\bsU^{r}_{g-r,\s}$.
Therefore $\sC$ is disjoint from the open substack
$\sU_{g-r}^{r,part}$. However, this contradicts the assumption
that $\pi(\sC)$ does not lie in $A_{g-r-1}^{Sat}$, and \DHrefpart{i} is proved.

For \DHrefpart{ii},
assume first that we are in positive characteristic
so that Keel's theorem \cite{Ke} is available.
Assume also, as an induction hypothesis, that
$L_h$ is semi-ample on $\sA_h^P$ for all $h<g$.

Suppose that 
$\cZ$ is an irreducible closed substack of $\sA_g^P$ with 
$\cZ.L_g^{\dim \cZ}=0$,
so that $\cZ$ lies in $\Exc(L_g)$. If $\sZ\cap \sA_g^{part}$
is not empty, then $\sZ$ lies in $\sA_{1,g-1}^P$,
as desired. So we can suppose
that $\sZ$ lies over $A_{g-r}^{Sat}$, but not over
$A_{g-r-1}^{Sat}$, where $r\ge 2$. Then $\sZ$ lies in some $\bsX_r$. 

Note first that $r<g$, since $L_g$ is ample on the fibre
$\pi^{-1}(A_0^{Sat})$ of $\pi$.
(This is the statement that $-D_g$ is $\pi$-ample, which holds
because, as a toroidal compactification, 
$\sA_g^P$ is defined by taking a convex hull.)

Once again we can use a torus action to construct a rational equivalence
$\sZ\sim\sZ_0=\sY+\sW$ where $\sY$ lies in $\bsU^r_{g-r,\s}$ and
the fibres of $\sW\to A_g^{Sat}$ are of strictly positive dimension.

Suppose first that $\sZ=\sY$, i.e., that $\sZ$ lies in some $\bsU^r_{g-r,\s}$.
Since $L_g$ is nef, we have $\sZ.L_g^{\dim \sZ}=0$. 
Put $\sZ^0 = \sZ\cap \sU_{g-r}^{r,part}$; this is open and dense in $\sZ$.

As in the proof of \DHrefpart{i},
it follows from
the linear equivalence \ref{(*)}
that the stable base locus of the restriction
$L_g\vert_{\bsU^r_{g-r,\s}}$ 
lies in the boundary $\bsU^r_{g-r,\s}\setminus\sU_{g-r}^{r,part}$.
Therefore, by
Kodaira's lemma (\cite {Ko} VI.2.15, VI.2.16),
the restriction $L'$ of $L_g$ to $\sU_{g-r}^{r,part}$ is semi-ample
but not big, so
that $\sZ^0$ is covered by open curves $\sC^0$ on which the morphism
defined by the linear system $\vert mL'\vert$, for some suitable 
integer $m$, is constant.

Then, for each $i$, the linear system
$\vert mL_{g-r+1}\vert$ defines a constant morphism on each
curve $\rho_i(\sC^0)$. By induction, $\vert mL_{g-r+1}\vert$
has no base points on $\sA_{g-r+1}^P$, so the closure of
$\rho_i(\sC^0)$ lies in $\Exc(L_{g-r+1})$.

By the induction
hypothesis, $\Exc(L_{g-r+1})=\sA_{1,g-r}^P$.
Moreover, if $\sU_{g-r}^{part}$ is identified with an open
substack of the boundary $D_{g-r+1}$ of $\sA_{g-r+1}^P$,
then taking the $j$-invariant of the elliptic factor to be $\infty$
shows that $\Exc(L_{g-r+1})\cap \sU_{g-r}^{part}$ 
contains the closure
of the zero-section of the semi-abelian scheme 
$\sU_{g-r}^{part}\to \sA_{g-r}^{part}$. 

Consider the intersection $\sI=\sA_{1,g-r}^P\cap \sU_{g-r}^{part}$,
taken inside $D_{g-r+1}$.

\begin{lemma} $\sI$ has just two irreducible components.
One is the locus of points of the form
$$(\infty, B, 0_B),$$
where $B\in \sA_{g-r}$ and $\infty$ is the point at infinity on $\sA_1^P$;
this is a copy of $\sA_{g-r}$. 
The other is the locus of points of the form
$$(E\times V, (0_E,v)),$$
where $E\in \sA_1$, $V\in \sA_{g-r-1}$ and $v\in V$ is arbitrary;
this is a copy of $\sA_1^P\times \sU_{g-r-1}$.
\begin{proof}
The only thing to notice is that on $\sA_{g-r+1}^{part}$, the
exceptional locus $\Exc(L_{g-r+1})$ includes 
the image of $\sA_1\times \sA_{g-r}^{part}$,
the locus where the cycle-theoretic image of the Kummer variety 
under the $2\theta$ linear system does not determine the abelian variety.
The pair $(V,v)$ corresponds to a compactification $\tV$ of 
some $\GG_m$-bundle over $V$, and then taking the Kummer variety
of $E\times \tV$ has the effect of ``losing the isomorphism class of $E$''.
\end{proof}
\end{lemma}

That is, for every $\txi_i:\sU^{r,part}_{g-r}\to\sU^{part}_{g-r}$, 
the image $\txi_i(\sZ)$ is contained in
the union of these two loci. Now consider the summand
$a\delta^*L_{g-r}$
that appears as a contribution 
to $L_g\vert_{\sU_{g-r}^{r,part}}$ in \ref{(*)}; since, by induction, $L_{g-r}$ is
semi-ample and $\Exc(L_{g-r})=\sA_{1,g-r-1}^P$,
this consideration shows that 
$$\sZ^0\cap \sU_{g-r}^r \subset \sA_1 \times \sU_{g-r-1}^r.$$
Now $\sZ$ lies in the closure of $\sZ^0\cap \sU_{g-r}^r$ in
$\bsU_{g-r,\s}^r$, so that
$\sZ$ is in (the image of)
$\sA_1^P \times \bsU_{g-r-1,\s}^r$ in $\sA_{1,g-1}^P$. In particular,
$\sZ$ lies in $\sA_{1,g-1}^P$.

Now drop the assumption that $\sZ$ lies in some $\bsU^r_{g-r,\s}$.
Then $\sZ$ specializes as above to $\sZ_0=\sW+\sY$ where $\sY$ lies in
some $\bsU^r_{g-r,\s}$. Since $L_g$ is nef, both $\sW$ and $\sY$
lie in $\Exc(L_g)$, so that, by what we have already proved, $\sY$ lies in $\sA_{1,g-1}^P$.

Recall that $\sZ$ lies in $\bsX_{r,\tau}$, the image of the closure of a $T$-bundle 
$\sT\to\sU_{g-r}^r$.
The specialization $\sZ\sim\sZ_0$ and the fact that $\sY$ lies in
$\sA_{1,g-1}^P$ show
that $\sZ$ lies in the image of the closure of the restriction of $\sT$ 
to the closed substack $\sA_1\times \sU_{g-r-1}^r$ of $\sU_{g-r}^r$.
But this restriction is of the form $\sA_1\times \sT_1$,
where $\sT_1$ is a $T$-bundle over $\sU_{g-r-1}^r$.
So $\sZ$ lies in 
$\sA_1^P\times \bsX'_{r,\tau}$, where $\bsX'_{r,\tau}$ is the image of
the closure
of $\sT_1$.
However, $\bsX'_{r,\tau}$ lies in $\sA_{g-1}^P$,
so that $\sZ$ lies in $\sA_{1,g-1}^P$,
as required.

That is, we have shown that in characteristic $p>0$,
$\Exc(L_g)\subset\sA_{1,g-1}^P$; the other inclusion is an
immediate consequence of the fact that $L_1$ is trivial. From Keel's theorem
we deduce that $L_g$ is semi-ample. It follows at once
that $\Exc(L_g)=\sA_{1,g-1}^P$ in characteristic zero,
and now \DHrefpart{ii}
is proved.

\DHrefpart{iii} 
We know now that $L_g$ is semi-ample, 
and it remains to show that the varieties $V_g$ behave as stated.

The multiplication morphism $\phi:\sA^P_{g-h}\times\sA^P_h\to \sA^P_g$
has the property that $\phi^*\sO(H_g)\cong\sO(H_{g-h})\boxtimes\sO(H_h)$
and $\phi^*\sO(D_g)\cong\sO(D_{g-h})\boxtimes\sO(D_h)$.
So  $\phi^*\sO(L_g)\cong\sO(L_{g-h})\boxtimes\sO(L_h)$.

Let $\psi_g:A_g^P\to V_g$ be the morphism defined by the linear system
$\vert mL_g\vert$ for sufficiently divisible $m$
and take $h=1$; since $L_1$ is trivial,
the composite
$$A_{g-1}^P\times A_1^P\to A^P_{g-1,1}\to A_g^P\to V_g$$
factors through $pr_1:A_{g-1}^P\times A_1\to A_{g-1}^P$, say
via $A_{g-1}^P\to V_g$. Let $A_{g-1}^P\to W_{g-1}$
denote the Stein factorization of this.
Since $A_{g-1}\to V_g$ is defined by some system $\vert mL_{g-1}\vert$,
it follows that $W_{g-1}$
is identified with $V_{g-1}$.
\end{proof}
\end{theorem}
\begin{corollary}\label{interior}
Suppose that $7\le g\le 10$ and that $\ch k=0$.
Then the first step in running the MMP on $A_g^P$
is the contraction $A_{g}^P\to V_{g}$.
Its fibres are the fibres of $A_{g-1,1}^P\to V_{g-1}$
and they meet the interior $A_g$ of $A_g^P$.
\begin{proof}
This follows at once from the results above
when it is recalled that $A_g^P$ has only terminal
singularities in this range [AS].
\end{proof}
\end{corollary}
\end{section}
\bibliography{alggeom,ekedahl}

\providecommand{\bysame}{\leavevmode\hbox to3em{\hrulefill}\thinspace}
\begin{thebibliography}{EGAIII:2}


\bibitem[AMRT]{AMRT} A.~Ash, D.~Mumford, M.~Rapoport and Y.-S.~Tai,
\emph{Smooth compactifications of locally symmetric varieties}, 2nd ed.,
CUP, 2010.

\bibitem[AS]{AS} J.~Armstrong and N.~Shepherd-Barron,
\emph{The singularities of $A_g^P$} arXiv:1604.05943

\bibitem[BC]{BC} E.S.~Barnes and M.J.~Cohn, \emph{On the inner product of positive
quadratic forms}, J. London Math. Soc. \textbf{12} (1975/76), 32--36.



\bibitem[DO]{DO}
I.~Dolgachev and D.~Ortland, \emph{Point sets in projective space and theta functions},
Ast\'erisque \textbf{165}, 1988.


\bibitem[FC]{FC}
G.~Faltings and C.~L. Chai, \emph{Degenerations of abelian varieties},
Springer Verlag, 1980.




\bibitem[H]{H}
R.W.H.T.~Hudson, \emph{Kummer's Quartic Surface}, 2nd ed., Cambridge, 1990.

\bibitem[HH]{HH}
B.~Hassett and D.~Hyeon, \emph{Log minimal model program for the
moduli space of stable curves: the first flip},
Annals of Math. \textbf{177} (2013), 911--968.



\bibitem[Ke]{Ke}
S.~Keel, \emph{Basepoint freeness for nef and big line bundles in positive 
characteristic}, Annals of Math., \textbf{149} (1999), 253--286.


\bibitem[KM]{KM}
S.~Keel and S.~Mori, \emph{Quotients by groupoids}, Annals of Math. \textbf{145}
(1997), 193--213.


\bibitem[Ko]{Ko}
J.~Koll{\'a}r, \emph{Rational curves on algebraic varieties}, Springer, 1996.


%


\bibitem[SB]{SB}
N.I.~Shepherd-Barron, \emph{Perfect forms and the moduli space
of abelian varieties}, Invent. Math. \textbf{163} (2006), 25--45.

\bibitem[S]{S} V.S.~Snurnikov, Ph.D. thesis, Cambridge, 2002.

\bibitem[T]{T} Y.-S.~Tai, \emph{On the Kodaira dimension of the moduli space 
of abelian varieties}, Invent. Math. \textbf{68} (1982), 425--439.


\end{thebibliography}
\bibliographystyle{pretex}
\end{document}